\documentclass[10 pt,reqno]{amsart}
\usepackage{amsfonts}
\usepackage{amsmath}
\usepackage{mathdots}
\usepackage{amsthm}
\usepackage{amssymb, cite}
\usepackage[backref=page]{hyperref}

\newcommand{\C}{\mathbb{C}}

\newcommand{\R}{\mathbb{R}}
\newcommand{\M}{\mathcal{M}}

\newcommand{\norm}[1]{\lVert #1 \rVert}

\newcommand{\bea}{\begin{eqnarray*}}
\newcommand{\eea}{\end{eqnarray*}}

\newcommand{\tr}{\operatorname{tr}}
\newcommand{\inner}[1]{\langle #1 \rangle}
\newcommand{\h}{\mathcal{H}}

\newcommand{\minimatrix}[4]{\begin{pmatrix} #1 & #2 \\ #3 & #4 \end{pmatrix}  }
\newcommand{\megamatrix}[9]{\begin{pmatrix} #1 & #2 & #3 \\ #4 & #5 & #6 \\ #7 & #8 & #9\end{pmatrix}  }

\allowdisplaybreaks

\newtheorem{Theorem}{Theorem}

\newtheorem{Proposition}[Theorem]{Proposition}
\newtheorem{Corollary}[Theorem]{Corollary}

\newtheorem{Lemma}[Theorem]{Lemma}
\newtheorem{Example}[Theorem]{Example}

\numberwithin{equation}{section}
\numberwithin{Theorem}{section}

\theoremstyle{definition}
\newtheorem*{Definition}{Definition}

\begin{document}

	\title[Unitary equivalence of a matrix to its transpose]
		{Unitary equivalence of a matrix to its transpose}

	\author{Stephan Ramon Garcia}
	\address{   Department of Mathematics\\
		Pomona College\\
		Claremont, California\\
		91711}
	\email{Stephan.Garcia@pomona.edu}
	\urladdr{http://pages.pomona.edu/\textasciitilde sg064747}

	\author{James E. Tener}
	\address{   Department of Mathematics\\
		University of California\\Berkeley, CA\\ 94721}
	\email{jtener@math.berkeley.edu}
	\urladdr{http://math.berkeley.edu/\textasciitilde jtener}

    
	\thanks{The first author was partially supported by NSF Grant DMS-0638789.  The second author
	was partially supported by a NSF Graduate Research Fellowship.}

	\begin{abstract}
		Motivated by a problem of Halmos, we obtain a canonical decomposition for complex matrices which
		are unitarily equivalent to their transpose (UET).  Surprisingly, the na\"ive assertion that 
		a matrix is UET if and only if it is unitarily equivalent to a complex symmetric matrix 
		holds for matrices $7 \times 7$ and smaller, but fails for matrices $8 \times 8$ and larger.
	\end{abstract}

\maketitle

\section{Introduction}

	In his problem book \cite[Pr.~159]{Halmos}, Halmos asks whether every square complex matrix is unitarily equivalent to
	its transpose (UET).  	For example, every finite Toeplitz matrix is unitarily equivalent to its transpose via the
	permutation matrix which reverses the order of the standard basis.
	Upon appealing to the Jordan canonical form, it follows that every square complex 
	matrix $T$ is \emph{similar} to its transpose $T^t$.  Thus
	similarity invariants are insufficient to handle Halmos' problem.  

	In his discussion, Halmos introduces the single counterexample
	\begin{equation}\label{eq-HalmosExample}
		\megamatrix{0}{1}{0}{0}{0}{2}{0}{0}{0}
	\end{equation}
	and proves that it is not UET.  A more recent example,
	due to George and Ikrahmov \cite{George}, is
	\begin{equation}\label{eq-George}
		\megamatrix{1}{0}{0}{4}{3}{0}{0}{2}{5}.
	\end{equation}
	While settling Halmos' original question in the negative,
	\eqref{eq-HalmosExample} and \eqref{eq-George} are only 
	isolated examples.  Our present aim is to obtain  
	a complete characterization and canonical decomposition
	of those matrices which are UET.
		
	From some perspectives, every square matrix is close to being UET.
	Indeed, it is known that for each $T \in M_n(\C)$ there exist unitary matrices $U$ and 
	$V$ such that $T = UT^tV$ \cite{Gordon}.  In fact, Taussky and Zassenhaus proved that
	for any field $\mathbb{K}$ each $T \in M_n(\mathbb{K})$ 
	is similar to its transpose \cite{Taussky} (for congruence see \cite{DDZ}).  	
	
	In addition to Halmos' question, this article is partially motivated by the recent explosion
	in work on linear preservers.
	In particular, linear maps on $M_n(\C)$ of the form $\phi(T) = UT^t U^*$ (where $U$ is unitary)
	feature prominently in the literature \cite{Chan,Clark,Gau, LPR, LR,LSN, Petek}.
	Our work completely classifies the fixed points of such maps.	

	Halmos' example \eqref{eq-HalmosExample} recently appeared in another context
	which also motivated the authors to consider his problem.	
	It is well-known that every square complex matrix is \emph{similar} to a 
	complex symmetric matrix (i.e., $T = T^t$) \cite[Thm.~4.4.9]{HJ}.  Indeed,
	various complex symmetric canonical forms have been proposed throughout the years
	\cite{Craven, Gantmacher, HJ, ScottCanonical, Wellstein}.	
	It turns out that Halmos' matrix \eqref{eq-HalmosExample} was among the first several
	matrices demonstrated to be \emph{not} unitarily equivalent to a complex
	symmetric matrix (UECSM) \cite[Ex.~1]{SNCSO} (see also \cite[Ex.~6]{CSOA} and \cite[Ex.~1]{Vermeer}).

	In \cite{UECSMGC} it was remarked in passing that
	\begin{equation}\label{eq-Implication}
		\text{$T$ is UECSM} \quad\Rightarrow\quad \text{$T$ is UET},
	\end{equation}
	raising the question of whether the converse is also true.  It turns out that 
	Vermeer had already studied this problem over $\R$ a few years earlier and provided an $8 \times 8$ 
	counterexample \cite[Ex.~3]{Vermeer}.  In fact, we prove that there are no smaller counterexamples:
	the implication \eqref{eq-Implication} can be reversed for matrices $7 \times 7$ and smaller,
	but not for matrices $8 \times 8$ or larger.	

	To state our main results, we require a few definitions.
	
	\begin{Definition}
		A $2d \times 2d$ block matrix of the form
		\begin{equation}\label{eq-SHM}
			T = \minimatrix{A}{B}{D}{A^t}
		\end{equation}
		where $B^t = -B$ and $D^t = - D$ is called
		\emph{skew-Hamiltonian} (SHM).
		A matrix $T \in M_{2d}(\C)$ that is unitarily equivalent to an
		skew-Hamiltonian matrix is called UESHM.
	\end{Definition}

	Needless to say, the matrices $A$, $B$, and $D$ in \eqref{eq-SHM} are necessarily
	$d \times d$.  A short computation reveals that if $T$ is skew-Hamiltonian, then
	\begin{equation}\label{eq-Symplectic}
		T = -\Omega T^t \Omega = \Omega T^t \Omega^*,
	\end{equation}
	where
	\begin{equation}\label{eq-Omega}
		\Omega = \minimatrix{0}{I}{-I}{0}.
	\end{equation}
	In particular, it follows immediately from \eqref{eq-Symplectic} that every SHM is UET.
	\medskip

	We are now ready to state our main result:

	\begin{Theorem}\label{TheoremMain}
		A matrix $T \in M_n(\C)$ is UET if and only if it is 
		unitarily equivalent to a direct sum of  (some of the summands may be absent):
		\begin{enumerate}\addtolength{\itemsep}{0.5\baselineskip}
			\item[I.] irreducible complex symmetric matrices (CSMs).
			\item[II.] irreducible skew-Hamiltonian matrices (SHMs).  Such matrices
				are necessarily $8 \times 8$ or larger.
			\item[III.] $2d \times 2d$ blocks of the form
				\begin{equation}\label{eq-TypeIII}
					\minimatrix{A}{0}{0}{A^t}
				\end{equation}
				where $A \in M_d(\C)$ is irreducible and neither UECSM nor UESHM.
				Such matrices are necessarily $6 \times 6$ or larger.
		\end{enumerate}
		Moreover, the unitary orbits of the three classes described above are pairwise disjoint.
	\end{Theorem}

	We use the term \emph{irreducible} in the operator-theoretic sense.
	Namely, a matrix $T \in M_n(\C)$ is called irreducible if $T$ is not unitarily
	equivalent to a direct sum $A \oplus B$ where $A \in M_d(\C)$ and $B \in M_{n-d}(\C)$
	for some $1 < d < n$.  Equivalently, $T$ is irreducible if and only if the only normal matrices
	commuting with $T$ are multiples of the identity.  In the following, we shall denote unitary
	equivalence by $\cong$.
	
\section{Some corollaries of Theorem \ref{TheoremMain}}
	
	The proof of Theorem \ref{TheoremMain} requires a number of preliminary results
	and is consequently deferred until Section \ref{SectionProof01}.   We focus here on
	a few immediate corollaries.

	\begin{Corollary}
		If $T \in M_n(\C)$ is irreducible and UET, then $T$ is either UECSM or UESHM.
	\end{Corollary}
	
	\begin{Corollary}
		If $n \leq 5$ and $T \in M_n(\C)$ is UET, then $T$ is unitarily equivalent to
		a direct sum of irreducible complex symmetric matrices.
	\end{Corollary}	
	
	Our next corollary implies that the converse of the implication \eqref{eq-Implication}
	holds for matrices $7 \times 7$ and smaller.  On the other hand,
	it is possible to show (see Section \ref{SectionIrreducible}) that the converse fails
	for matrices $8 \times 8$ and larger.

	\begin{Corollary}\label{Corollary7x7}
		If $n \leq 7$ and $T \in M_n(\C)$ is UET, then $T$ is UECSM.
	\end{Corollary}

	\begin{proof}
		By Theorem \ref{TheoremMain}, $T$ is unitarily equivalent to a direct sum of blocks of
		type I or III.   It turns out (see Lemma \ref{LemmaBlock}) that any matrix of the form \eqref{eq-TypeIII}
		is UECSM (although clearly not irreducible and hence not of Type I) whence $T$ is
		itself UECSM. 
	\end{proof}

	We close this section with a few remarks about $2 \times 2$ and $3 \times 3$
	matrices.  For $2 \times 2$ matrices it has long been known \cite{Murnaghan} that
	$A \cong B$ if and only if $\Phi(A) = \Phi(B)$ where $\Phi:M_2(\C)\to \C^3$ is the function
	\begin{equation*}
		\Phi(X) = (\tr X, \tr X^2, \tr X^*X).
	\end{equation*}
	Since $\Phi(X) = \Phi(X^t)$ for all $X \in M_2(\C)$ we immediately obtain
	the following useful lemma from Corollary \ref{Corollary7x7}.
	
	\begin{Lemma}\label{Lemma2x2UECSM}
		Every $2 \times 2$ matrix is UECSM.
	\end{Lemma}

	For other proofs of the preceding lemma see
	\cite[Cor.~3]{UECSMGC}, \cite[Cor.~3.3]{Chevrot}, \cite{UECSMMC}, \cite[Ex.~6]{CSOA}, \cite[Cor.~1]{SNCSO}, \cite[p.~477]{McIntosh}, or
	\cite[Cor.~3]{Tener}.

	Based on Specht's Criterion \cite{Specht}, Pearcy obtained a list 
	of nine words in $X$ and $X^*$ so that $A,B \in M_3(\C)$
	are unitarily equivalent if the traces of these words are equal
	for $X = A$ and $X = B$ \cite{Pearcy}.  Later Sibirski{\u\i} \cite{Sibirskii}
	showed that two of these words are unnecessary and, moreover, 
	that $A \cong B$ if and only if $\Phi(A) = \Phi(B)$ 
	where $\Phi:M_3(\C)\to \C^7$ is the function defined by
	\begin{equation}\label{eq-f(X)}
		\Phi(X) = (\tr X, \, \tr  X^2,\, \tr X^3,\, \tr X^* X,\, \tr X^*X^2, \, \tr {X^*}^2 X^2, \, \tr X^* X^2 {X^*}^2X).
	\end{equation}
	The preceding can be used to develop a simple test for checking whether a $3 \times 3$
	matrix is UECSM.  Some general approaches to this problem in higher dimensions
	can be found in	\cite{UECSMGC, UECSMMC, Tener} (see also \cite[Thm.~3]{Vermeer}).
	 
	\begin{Proposition}\label{Proposition3x3}
		If $X \in M_3(\C)$, then $X$ is UECSM if and only if 
		\begin{equation}\label{eq-Trace3}
			\tr X^* X^2 {X^*}^2 X = \tr X {X^*}^2 X^2 X^*.
		\end{equation}
	\end{Proposition}

	\begin{proof}
		Since a $3 \times 3$ matrix is UECSM if and only if it is UET (by Theorem \ref{TheoremMain}),
		it suffices to prove that \eqref{eq-Trace3} is equivalent to asserting
		that $X \cong X^t$.  This in turn is equivalent to proving that
		$\Phi(X) = \Phi(X^t)$ for the function $\Phi:M_3(\C) \to \C^7$ defined by \eqref{eq-f(X)}.
		Let $\phi_1(X), \phi_2(X), \ldots, \phi_7(X)$ denote the entries of \eqref{eq-f(X)} and note
		that $\phi_i(X) = \phi_i(X^t)$ for $i = 1,2,3$.  Moreover, a short computation shows that
		\begin{align*}
			\phi_4(X) &= \tr X^*X = \tr XX^* = \tr {X^*}^t X^t = \tr {X^t}^*X^t = \phi_4(X^t), \\
			\phi_5(X) &= \tr X^*X^2 = \tr X^2 X^* = \tr {X^*}^t {X^2}^t = \tr {X^t}^* {X^t}^2 = \phi_5(X^t), \\
			\phi_6(X) &= \tr {X^*}^2 X^2 = \tr (X^2)^t {(X^*)^2}^t = \tr (X^t)^2 {(X^t)^*}^2
			=  \tr {(X^t)^*}^2 (X^t)^2 = \phi_6(X^t).
		\end{align*}
		Thus $X$ is UECSM if and only if $\phi_7(X) = \phi_7(X^t)$.  Since
		\begin{equation*}
		\phi_7(X^t) = \tr {X^t}^* {X^t}^2 {(X^t)^*}^2 X^t = X {X^*}^2 X^2 X^*
		\end{equation*}
		the desired result follows.
	\end{proof}

	Pearcy also proved that one need only check traces of words of
	length $\leq 2n^2$ to test two $n \times n$ matrices for unitary equivalence \cite[Thm.~2]{Pearcy2}.
	As Corollary \ref{Corollary7x7} and Proposition \ref{Proposition3x3} suggest, 
	a similar algorithm could be developed for $n \leq 7$.
	However, the number of words one must consider grows too rapidly for this
	approach to be practical even in the $4 \times 4$ case.

\section{Building blocks of type I: UECSMs}\label{SectionUECSM}
	In this section we gather together some information about UECSMs
	that will be necessary in what follows.  We first require the notion of a conjugation.
	
	\begin{Definition} 
		A function $C: \C^n\to\C^n$ is called a
		\emph{conjugation} if it has the following three properties:
		\begin{enumerate}\addtolength{\itemsep}{0.5\baselineskip}
			\item $C$ is conjugate-linear,
			\item $C$ is isometric:  $\inner{ x,y} = \inner{Cy,Cx}$ for all $x,y \in \C^n$,
			\item $C$ is involutive:  $C^2 = I$.
		\end{enumerate}
	\end{Definition}
	
	In light of the polarization identity, condition (ii) in the preceding definition
	is equivalent to asserting that $\norm{Cx} = \norm{x}$ for all $x \in \C^n$.	
	Let us now observe that $T$ is a complex symmetric matrix (i.e., $T = T^t$) if and only if 
	$T = JT^*J$, where $J$ denotes the \emph{canonical conjugation}
	\begin{equation}\label{eq-Canonical}
		J( z_1,z_2, \ldots, z_n) = ( \overline{z_1} , \overline{z_2}, \ldots, \overline{z_n})
	\end{equation}
	on $\C^n$.  Moreover, we also have
	\begin{align}
		\overline{T} &= JTJ, \label{eq-JTJ} \\
		T^t &= JT^*J,\label{eq-JT*J}
	\end{align}
	where $\overline{T}$ is the entry-by-entry complex conjugate of $T$ (these relations have long
	been used in the study of Hankel operators \cite{Butz}).

	\begin{Lemma}\label{LemmaS}
		$C$ is a conjugation on $\C^n$ if and only if $C = UJ$ where $U$ is a complex
		symmetric (i.e., $U = U^t$) unitary matrix.
	\end{Lemma}

	\begin{proof}
		If $C$ is a conjugation on $\C^n$, then $U = CJ$ is an isometric linear
		map and hence unitary.  It follows from \eqref{eq-JTJ} 
		that $U\overline{U} = UJUJ = C^2 = I$ whence $\overline{U} = U^*$ so that
		$U = U^t$ as claimed.  Conversely, if $U$ is a complex symmetric unitary matrix then
		$C = UJ$ is conjugate-linear, isometric, and satisfies $C^2 = I$ by a similar computation.
	\end{proof}
	
	The relevance of conjugations to our work lies in the following lemma (the equivalence of
	(i) and (ii) was first noted in \cite[Thm.~3]{Vermeer}).

	\begin{Lemma}\label{LemmaUECSM}
		For $T \in M_n(\C)$, the following are equivalent:
		\begin{enumerate}\addtolength{\itemsep}{0.5\baselineskip}
			\item $T$ is UECSM,
			\item $T = UT^t U^*$ for some complex symmetric unitary matrix $U$,
			\item $T = CT^*C$ for some conjugation $C$ on $\C^n$.
		\end{enumerate}
		\smallskip
		In particular, if $T$ is UECSM, then $T$ is UET.  
	\end{Lemma}

	\begin{proof}
		(i) $\Rightarrow$ (ii) If $Q$ is unitary and $Q^*TQ = S$ is complex symmetric, then
		$Q^*TQ = (Q^*TQ)^t = Q^t T^t \overline{Q}$ whence
		$T = UT^t U^*$ where $U = QQ^t$.
		\medskip

		\noindent (ii) $\Rightarrow$ (iii)
		If $T = UT^t U^*$ where $U = U^t$ is unitary,
		then $C = UJ$ is a conjugation by Lemma \ref{LemmaS}.  Since
		$U = CJ$ and $U^* = JC$, it follows from \eqref{eq-JT*J} that
		$T = UT^t U^* = CJT^tJC = CT^*C$.
		\medskip

		\noindent (iii) $\Rightarrow$ (i)
		Suppose that $T = CT^*C$ for some conjugation $C$ on $\C^n$.
		By \cite[Lem.~1]{CSOA} there exists 
		an orthonormal basis $e_1,e_2,\ldots,e_n$ such that $Ce_i = e_i$
		for $i = 1,2,\ldots,n$.  Let $Q = (e_1 | e_2| \cdots |e_n)$ be the
		unitary matrix whose columns are these basis vectors.  The
		matrix $S = Q^*TQ$ is complex symmetric since the $ij$th 
		entry $[S]_{ij}$ of $S$ satisfies
		\begin{equation*}
			[S]_{ij} = \inner{T e_j,e_i} = \inner{CT^*Ce_j,e_i} 
			= \inner{e_i, T^*e_j} = \inner{Te_i,e_j} = [S]_{ji}. \qedhere
		\end{equation*}
	\end{proof}
	
	In order to obtain the decomposition guaranteed by
	Theorem \ref{TheoremMain}, we must be able to break apart matrices which are
	UECSM into direct sums of simpler matrices.  Unfortunately,
	the class UECSM is not closed under restrictions to direct summands
	(see Example \ref{Example6x6} in Section \ref{SectionAAt}),
	making our task more difficult.  We begin by considering a special case 
	where such a reduction is possible.
	
	\begin{Lemma}\label{LemmaUECSMReduce01}
		Suppose that $T \in M_n(\C)$ is UECSM and let $C$ be a conjugation
		for which $T = CT^*C$.  If $\M$ is a proper, nontrivial subspace of $\C^n$ 
		that reduces $T$ and satisfies $\M \cap C\M \neq \{ 0\}$, then 		
		$T \cong T_1 \oplus T_2$ where $T_1$ and $T_2$ are UECSM.
	\end{Lemma}
	
	\begin{proof}		
		Let us initially regard $T$ as a linear operator $T:\C^n \to \C^n$
		and work in a coordinate-free manner until the conclusion of the proof.
		Since $\M$ reduces $T = CT^*C$, a short computation reveals that
		\begin{equation}\label{eq-TCM}
			T(C\M) = (TC)\M = (CT^*)\M = C(T^*\M) \subseteq C\M
		\end{equation}
		whence $C\M$ is $T$-invariant.  Upon replacing $T$ with $T^*$ in \eqref{eq-TCM}
		we find that $C \M$ is also $T^*$-invariant whence $C\M$ reduces $T$.  It follows
		that the subspaces
		\begin{equation*}
			\h_1 = \M\cap C\M , \qquad \h_2 = \h_1^{\perp},
		\end{equation*}
		both reduce $T$.   Observe that $\h_1,\h_2$ are both proper and nontrivial by
		assumption.  Moreover, note that $\h_1$ is $C$-invariant by definition and, since $C^2 = I$,
		it follows that $C\h_1 = \h_1$.  In light of the fact that $C$ is isometric it also follows
		that $C\h_2 = \h_2$.  With respect to the orthogonal decomposition 
		$\C^n = \h_1\oplus \h_2$, we have the block-operator decompositions
		$T = T_1 \oplus T_2$ and $C = C_1 \oplus C_2$ where $C_1,C_2$ are
		conjugations on $\h_1, \h_2$, respectively.  Expanding
		the identity $T = CT^*C$ in terms of block operators reveals that
		$T_1 = C_1 T_1^*C_1$ and $T_2 = C_2 T_2^* C_2$.
		Upon identifying $T_1,T_2$ as matrices computed with respect to 
		some orthonormal bases of $\h_1,\h_2$, respectively, 
		we conclude that $T_1$ and $T_2$ are UECSM by Lemma \ref{LemmaUECSM}.
	\end{proof}

	If a matrix which is UECSM is reducible and the preceding lemma does not apply, then we have the following decomposition.

	\begin{Lemma}\label{LemmaUECSMReduce02}
		Suppose that $T = A \oplus B \in M_n(\C)$ is UECSM, $A \in M_d(\C)$, 
		$B \in M_{n-d}(\C)$, and let $C$ be a conjugation on $\C^n$
		for which $T = CT^*C$.  If $\M \cap C\M = \{0\}$
		for every proper reducing subspace $\M$ of $T$, then $n = 2d$ 
		and $T \cong A \oplus A^t$ where $A \in M_d(\C)$ is irreducible.
	\end{Lemma}
	
	\begin{proof}
		Since $T = A \oplus B$, the subspaces
		$\M = \C^d \oplus \{0\}$ and $\M^{\perp} = \{0\} \oplus \C^{n-d}$ reduce $T$.
		By hypothesis and the fact that $C(\M^{\perp}) = (C\M)^{\perp}$ it follows that
		\begin{align}
			\M \cap C\M &= \{0\}, \label{eq-MM01} \\
			\M^{\perp} \cap C(\M^{\perp}) &= \{0\}. \label{eq-MM02}
		\end{align}
		Since $\dim C\M = \dim \M = d$, it follows from \eqref{eq-MM01} that
		$2d \leq n$.  Similarly, since $\dim C(\M^{\perp}) = \dim \M^{\perp} = n-d$, it follows
		from \eqref{eq-MM02} that $2(n-d) \leq n$.  Putting these two inequalities together
		reveals that $n = 2d$. In fact, a similar argument shows that \emph{every} proper,
		nontrivial reducing subspace of $T$ is of dimension $d$.
		Moreover, it also follows that $A$ and $B$ are both irreducible, since otherwise $T$
		would have a nontrivial reducing subspace of dimension $< d$.

		Letting
		\begin{equation}\label{eq-PIO}
			P = \minimatrix{I}{0}{0}{0}
		\end{equation}
		denote the orthogonal projection onto $\M$, we note that $R = CPC$ is the orthogonal projection
		onto $C\M$.  Moreover, since $T = CT^*C$ it follows that $C\M$ reduces $T$.
		Writing
		\begin{equation}\label{eq-RProj}
			R = \minimatrix{R_{11} }{R_{12}}{R_{12}^*}{ R_{22} },
		\end{equation}
		where $R_{11}^*=R_{11}$ and $R_{22}^*=R_{22}$, and then expanding out 
		the equation $RT = TR$ block-by-block we find that
		$AR_{11} = R_{11}A$, $BR_{22} = R_{22}B$, and 
		\begin{equation}\label{eq-ARRB}
			AR_{12} = R_{12}B.
		\end{equation}
		Since $A$ and $B$ are both irreducible, it follows that $R_{11} = \alpha I$ and $R_{22} = \beta I$ where
		$0 \leq \alpha , \beta \leq 1$ (since $R$ is an orthogonal projection).
		Since $R^2  = R$ we also have
		\begin{equation}\label{eq-RRR}
			\minimatrix{ \alpha^2 I + R_{12} R_{12}^* }{ (\alpha + \beta) R_{12} }
			{ (\alpha + \beta) R_{12}^*}{ \beta^2I + R_{12}^*R_{12} }
			= \minimatrix{\alpha I }{ R_{12} }{ R_{12}^*}{\beta I},
		\end{equation}
		which presents three distinct possibilities.
		
		\medskip
		\noindent\textbf{Case 1}:  Suppose that $\alpha = 1$.  
		Looking at the $(1,1)$ entry in \eqref{eq-RRR}, we find that $R_{12} R_{12}^* = 0$
		whence $R_{12} = 0$.  From the $(2,2)$ entry in \eqref{eq-RRR} we find that $\beta^2 = \beta$ whence
		either $\beta = 0$ or $\beta = 1$.  Both cases are easy to dispatch.
		\begin{itemize}\addtolength{\itemsep}{0.5\baselineskip}
			\item If $\beta = 0$, then $R = P$,
				the orthogonal projection \eqref{eq-PIO} onto $\M = \C^d \oplus \{0\}$.
				Since $R$ is the orthogonal projection onto $C\M$ we have
				$\M = C\M$, which contradicts the hypothesis that $\M \cap C \M = \{0\}$.

			\item If $\beta = 1$, then $R = I$ which contradicts the fact that $\dim C\M = d$.
		\end{itemize}
  
		\medskip
		\noindent\textbf{Case 2}:  Suppose that $\alpha = 0$.  As before, we find that 
		$R_{12} = 0$ and $\beta^2 =\beta$.  Since $\beta = 0$ leads to the contradiction
		$R = 0$, it follows that $\beta = 1$ and
		\begin{equation*}
			R = \minimatrix{0}{0}{0}{I},
		\end{equation*}
		the orthogonal projection onto $\M^{\perp} = \{0\} \oplus \C^d$ (i.e., $C \M = \M^{\perp}$).
		With respect to the orthogonal decomposition $\C^{2n} = \C^d \oplus \C^d$ we have
		\begin{equation*}
			C = \minimatrix{0}{UJ}{U^tJ}{0}
		\end{equation*}
		where $U$ is a unitary matrix and $J$ is the canonical conjugation
		on $\C^d$ (by Lemma \ref{LemmaS}).  Expanding the equality $T = CT^*C$
		block-by-block reveals that $A = UJB^*U^tJ  = UB^tU^*$,
		from which we conclude that $A \cong B^t$.  Thus $T \cong A \oplus A^t$, as claimed.

		\medskip
		\noindent\textbf{Case 3}:  Suppose that $0 < \alpha < 1$.  In this case
		an examination of the $(1,1)$ entry in \eqref{eq-RRR} reveals that $R_{12} \neq 0$.  Looking next
		at the $(1,2)$ entry in \eqref{eq-RRR} we find that $\beta = 1 - \alpha$ from which 
		\begin{equation*}
			R_{12}^* R_{12} = R_{12} R_{12}^* = \alpha( 1 - \alpha) I
		\end{equation*}
		follows upon consideration of the $(1,1)$ and $(2,2)$ entries of \eqref{eq-RRR}.
		In other words, $R_{12} = \sqrt{ \alpha (1 - \alpha) } U^*$ for some $d \times d$ unitary matrix $U$ so that
		\begin{equation*}
			R = \minimatrix{ \alpha I}{ \sqrt{ \alpha(1 - \alpha) } U^* }{ \sqrt{ \alpha(1 - \alpha) } U}{ 1 - \alpha}.
		\end{equation*}
		By \eqref{eq-ARRB} it also follows that $UA = BU$ whence $A \cong B$.

		Now recall that $R=CPC$ is the orthogonal projection onto $C\M$ and that $C = SJ$ for some
		$n \times n$ complex symmetric unitary matrix
		\begin{equation*}
			S = \minimatrix{S_{11} }{ S_{12} }{ S_{12}^t }{ S_{22} }.
		\end{equation*}
		Writing $CP = RC$ as $SJP = RSJ$ (where $J$ denotes the canonical conjugation
		on $\C^n$) we note that $JP = PJ$
		by \eqref{eq-PIO} and conclude that $SP = RS$.  
		In other words,
		\begin{equation}\label{eq-BigMess}\small
			\minimatrix{ S_{11} }{0}{ S_{12}^t }{0} = 
			\minimatrix{ \alpha S_{11} + \sqrt{ \alpha( 1 - \alpha) } U^* S_{12}^t }
			{\alpha S_{12} + \sqrt{ \alpha(1 - \alpha) } U^* S_{22} }
			{\sqrt{ \alpha(1 - \alpha) } U S_{11} + (1 - \alpha) S_{12}^t}
			{\sqrt{\alpha (1 - \alpha)} U S_{12} + (1 - \alpha) S_{22}}.
		\end{equation}
		Examining the $(1,1)$ entry of \eqref{eq-BigMess} and using the fact that 
		$S_{11}^t = S_{11}$ we find that
		\begin{equation}\label{eq-S121}
			S_{12} = \sqrt{ \tfrac{ 1 - \alpha}{ \alpha} } S_{11} U^t.
		\end{equation}
		The $(1,2)$ entry of \eqref{eq-BigMess} now tells us that
		\begin{equation}\label{eq-S122}
			S_{12} = - \sqrt{ \tfrac{ 1 - \alpha}{\alpha} } U^* S_{22}.
		\end{equation}
		From \eqref{eq-S121} and \eqref{eq-S122} we have $S_{22} = - US_{11} U^t$ and hence
		\begin{equation*}
			S = \minimatrix{ S_{11} }{  \sqrt{ \tfrac{1-\alpha}{\alpha} } S_{11} U^t}
			{  \sqrt{ \tfrac{1-\alpha}{\alpha} } U S_{11}}{- US_{11}U^t}.
		\end{equation*}
		Recalling that $S = S^t$ is unitary, we have $I = S^*S = SS^*$,
		which can be expanded block-by-block to reveal that $S_{11} S_{11}^* = S_{11}^* S_{11} = \alpha I$.
		In other words, $S_{11} = \sqrt{ \alpha} \,V$ for some $d \times d$ complex symmetric unitary matrix $V$.
		Thus $S$ assumes the form
		\begin{equation*}
			S = \minimatrix{ \sqrt{\alpha}\, V }{  \sqrt{ 1-\alpha }\, V U^t}
			{  \sqrt{ 1-\alpha }\, UV}{- \sqrt{\alpha}\,UVU^t}.
		\end{equation*}
		Since $TC = CT^*$ and $C = SJ$ we have $TS = ST^t$, which yields
		\begin{equation*}
			\minimatrix{ \sqrt{\alpha} \,AV}{ \sqrt{1- \alpha}\,AVU^t}{\sqrt{1 - \alpha}\,BUV}{ -\sqrt{\alpha} BUVU^t}
			= \minimatrix{\sqrt{\alpha}\,VA^t}{ \sqrt{1 - \alpha}\,VU^tB^t}{ \sqrt{1 - \alpha}\,UVA^t}{-\sqrt{\alpha}\,UVU^tB^t}.
		\end{equation*}
		From the preceding we find $A \cong A^t$, $B \cong B^t$ (i.e., $A$ and $B$ are UET),
		and $A \cong B^t$.  In particular, $T \cong A \oplus A^t$, as claimed.
	\end{proof}
	
	Putting the preceding lemmas together we obtain the following proposition.
	
	\begin{Proposition}\label{PropositionUECSM}
		If $T \in M_n(\C)$ is UECSM, then $T$ is unitarily equivalent to a
		direct sum of matrices, each of which is either
		\begin{enumerate}\addtolength{\itemsep}{0.5\baselineskip}
			\item an irreducible complex symmetric matrix,
			\item a block matrix of the form $A \oplus A^t$ where $A$ is irreducible.
		\end{enumerate}
	\end{Proposition}
	
	\begin{proof}
		We proceed by induction on $n$.  The case $n=1$ is trivial.
		For our induction hypothesis, suppose that the theorem
		holds for $M_k(\C)$ with $k = 1,2,\ldots,n-1$.  Now suppose that $T \in M_n(\C)$
		is UECSM.  
		If $T$ is irreducible, then it is already of the form (i) and there is
		nothing to prove.  Let us therefore assume that $T$ is reducible.
		There are now two possibilities:
		\medskip
		
		\noindent \textsc{Case 1}:
		If $T$ has a proper, nontrivial reducing subspace $\M$ such that $\M \cap C\M \neq \{0\}$, 
		then Lemma \ref{LemmaUECSMReduce01} asserts that $T \cong A \oplus B$ 
		where $A$ and $B$ are UECSM.
		The desired conclusion now follows by the induction hypothesis.
		\medskip
		
		\noindent\textsc{Case 2}:		
		If every proper, nontrivial reducing subspace $\M$ of $T$
		satisfies $\M \cap C\M = \{0\}$, then $T$ satisfies the hypotheses of 
		Lemma \ref{LemmaUECSMReduce02}.  Therefore $n =2d$ and 
		$T \cong A \oplus A^t$ for some irreducible $A \in M_d(\C)$.
	\end{proof}

\section{Building block II:  UESHMs}\label{SectionSHM}
	
	As we noted in the introduction, there exist matrices which are UET but \emph{not} UECSM.
	In order to characterize those matrices which are UET, 
	we must introduce a new family of matrices along with the following definition.

	\begin{Definition} 
		A function $K: \C^n\to\C^n$ is called an
		\emph{anticonjugation} if it satisfies the following three properties:
		\begin{enumerate}\addtolength{\itemsep}{0.5\baselineskip}
			\item $K$ is conjugate-linear,
			\item $K$ is isometric:  $\inner{ x,y} = \inner{Ky,Kx}$ for all $x,y \in \C^n$,
			\item $K$ is skew-involutive:  $K^2 = -I$.
		\end{enumerate}
	\end{Definition}
	
	Henceforth, the capital letter $K$ will be reserved exclusively to denote anticonjugations.
	The proof of the following lemma is virtually identical to that of Lemma \ref{LemmaS} 
	and is therefore omitted.

	\begin{Lemma}\label{LemmaSJ}
		$K$ is an anticonjugation on $\C^n$ if and only if $K = SJ$ where $S$ is a 
		skew-symmetric (i.e., $S = -S^t$) unitary matrix.
	\end{Lemma}

	Unlike conjugations, which can be defined on spaces of arbitrary dimension,
	anticonjugations can act only on spaces of \emph{even} dimension.

	\begin{Lemma}\label{LemmaOdd}
		If $n$ is odd, then there does not exist an anticonjugation $K$ on $\C^n$.
	\end{Lemma}

	\begin{proof}
		By Lemma \ref{LemmaSJ} it suffices to prove that there 
		does not exist an $n \times n$ skew-symmetric unitary matrix $S$.
		If $S$ is such a matrix, then $\det S = \det (S^t) = \det(-S) = (-1)^n \det S = -\det S$	
		since $n$ is odd.  This implies that $\det S = 0$, a contradiction.
	\end{proof}
	
	Observe that skew-symmetric unitaries and their corresponding 
	anticonjugations exist when $n = 2d$ is even.  For
	example, let $S$ in Lemma \ref{LemmaSJ} be given by
	\begin{equation*}
		S = \bigoplus_{i=1}^{d} \minimatrix{0}{1}{-1}{0}.
	\end{equation*}

	\begin{Lemma}\label{LemmaONB}
		If $K$ is an anticonjugation on $\C^{2d}$, then
		there is an orthonormal basis $e_1,e_2,\ldots, e_{2d}$ of $\C^{2d}$ such that
		\begin{equation}\label{eq-CanAnti}
			Ke_i = 
			\begin{cases}
				e_{i+d} & \text{if $1 \leq i \leq d$}, \\
				-e_{i-d} & \text{if $d+1 \leq i \leq 2d$}.
			\end{cases}
		\end{equation}
	\end{Lemma}

	\begin{proof}
		The desired basis can be constructed inductively.
		First note that $\inner{x,Kx} = \inner{K^2x ,Kx} = - \inner{x,Kx}$ whence
		\begin{equation}\label{eq-xKx}
			\inner{x,Kx} = 0	
		\end{equation}
		for every $x \in \C^{2d}$.
		Now let $e_1$ be any unit vector, set $e_{d+1} = Ke_1$, and note that
		$Ke_{d+1} = -e_1$ since $K^2 = -I$.  In light of \eqref{eq-xKx}
		the set $\{ e_1, e_{d+1} \}$ is orthonormal.  Next select a unit vector $e_2$ in
		$\{e_1,e_{d+1}\}^{\perp}$ and let $e_{d+2} = Ke_2$.  A few additional computations
		reveal that the set $\{e_1,e_2,e_{d+1},e_{d+2}\}$ is orthonormal:
		\begin{align*}
			\inner{e_{d+2},e_1} &= \inner{Ke_1, Ke_{d+2}} = -\inner{e_{d+1}, e_2} = 0, \\
			\inner{e_{d+1},e_{d+2} } &= \inner{ Ke_{d+2}, Ke_{d+1} } = \inner{e_2,e_1} = 0.
		\end{align*}
		Continuing in this fashion we obtain the desired orthonormal basis.		
	\end{proof}

	\begin{Lemma}\label{LemmaUESHM}
		For $T \in M_n(\C)$ the following are equivalent:
		\begin{enumerate}\addtolength{\itemsep}{0.5\baselineskip}
			\item $T$ is UESHM,
			\item $T = UT^tU^*$ for some skew-symmetric unitary matrix $U$,
			\item $T = -KT^*K$ for some anticonjugation $K$ on $C^n$.
		\end{enumerate}		
		\smallskip
		In particular, if $T$ is UESHM, then $T$ is UET.  
		Moreover, for any of (i), (ii), or (iii) to hold $n$ must be even.
	\end{Lemma}

	\begin{proof}
		(i) $\Rightarrow$ (ii) Suppose that $Q^*TQ = S$ where $Q$ is unitary and $S$
		is skew-Hamiltonian.  By \eqref{eq-Symplectic} we have $S = \Omega S^t \Omega^t$ where
		$\Omega$ denotes the matrix \eqref{eq-Omega}.  It follows that
		\begin{equation*}
			Q^*TQ = S = \Omega S^t \Omega^t = (\Omega S \Omega^t)^t
			= (\Omega Q^*TQ \Omega^t)^t = \Omega Q^t T^t \overline{Q} \Omega^t
		\end{equation*}
		whence $T = UT^tU^*$ where $U = Q\Omega Q^t$ is a skew-symmetric unitary matrix.
		\medskip

		\noindent (ii) $\Rightarrow$ (iii)
		Suppose that $T = UT^t U^*$ where $U = -U^t$ is a unitary matrix.
		We claim that $K = UJ$ is an anticonjugation.
		Indeed, $K$ is conjugate-linear, isometric, and satisfies
		\begin{equation*}
			K^2 = (UJ)(UJ) = U(JUJ) = U\overline{U} = U(-U^*) = - I.
		\end{equation*}
		Putting this all together we find that
		\begin{equation*}
			T = UT^t U^* = U(JT^*J)(-\overline{U})
			= - (UJ)T^*(UJ) = - KT^*K,
		\end{equation*}
		as desired.
		\medskip

		\noindent (iii) $\Rightarrow$ (i)
		By Lemma \ref{LemmaOdd} we know that $n$ is even, say $n = 2d$.  Let
		$e_1,e_2,\ldots,e_{2d}$ be the orthonormal basis provided by Lemma \ref{LemmaONB}
		and let $Q = (e_1 | e_2| \cdots | e_{2d})$ be the $2d\times 2d$ unitary matrix whose columns
		are the basis vectors $e_1,e_2,\ldots,e_{2d}$.  The $ij$th entry $[S]_{ij}$ 
		of the matrix $S = Q^*TQ$ is given by the formula
		\begin{align*}
			[S]_{ij}
			&= \inner{Te_j,e_i} 
			= \inner{ -KT^*Ke_j,e_i} 
			= \inner{Ke_i,T^*Ke_j}
			= \inner{TKe_i,Ke_j} \\
			&=
			\left\{
			\begin{array}{rcll}
				\inner{Te_{i+d},e_{j+d} } &=& [S]_{j+d,i+d}  & \text{if $1 \leq i,j \leq d$}, \\[3pt]
				-\inner{Te_{i+d},e_{j-d} } &=& -[S]_{j-d,i+d} & \text{if $1 \leq i \leq d < j \leq 2d$}, \\[3pt]
				-\inner{Te_{i-d},e_{j+d} } &=& -[S]_{j+d,i-d} & \text{if $1 \leq j \leq d < i \leq 2d$}, \\[3pt]
				\inner{Te_{i-d},e_{j-d} } &=& [S]_{j-d,i-d} & \text{if $d \leq i,j \leq 2d$}.
			\end{array}
			\right.
		\end{align*}
		In other words, $S$ is of the form \eqref{eq-SHM} whence $T$ is UESHM.
	\end{proof}
	
	Putting the preceding material together, we obtain the following:

	\begin{Proposition}\label{PropositionUESHM}
		If $T \in M_n(\C)$ is UESHM, then $T$ is unitarily equivalent to a
		direct sum of matrices, each of which is either
		\begin{enumerate}\addtolength{\itemsep}{0.5\baselineskip}
			\item an irreducible skew-Hamiltonian matrix,
			\item a block matrix of the form $A \oplus A^t$ where $A$ is irreducible.
		\end{enumerate}
	\end{Proposition}
	
	\begin{proof}
		The proofs of Lemmas \ref{LemmaUECSMReduce01} and \ref{LemmaUECSMReduce02}
		go through \emph{mutatis mutandis} with $K$ in place of $C$ and $T = - KT^*K$ in place of
		$T = CT^*C$.  One then mimics the proof of Proposition \ref{PropositionUECSM}.
	\end{proof}

Before proceeding, let us remark that it is possible to show that every skew-Hamiltonian matrix
is \emph{similar} to a block matrix of the form $A \oplus A^t$ \cite{Waterhouse}.

\section{Building block III: matrices of the form $A \oplus A^t$}\label{SectionAAt}

	One of the difficulties in establishing Theorem \ref{TheoremMain} is the fact that 
	there is a nontrivial overlap between the classes UECSM and UESHM.	
	As Propositions \ref{PropositionUECSM} and \ref{PropositionUESHM} suggest,
	matrices of the form $A \oplus A^t$ are both UECSM and UESHM.

	\begin{Lemma}\label{LemmaBlock}
		If $A \in M_d(\C)$, then 
		\begin{equation}\label{eq-UET2x2}
			T = \minimatrix{A}{0}{0}{A^t} 
		\end{equation}	
		is both UECSM and UESHM.  In particular, 
		any matrix of the form \eqref{eq-UET2x2} is UET.
	\end{Lemma}

	\begin{proof}
		Letting $J$ denote the canonical conjugation \eqref{eq-Canonical} on $\C^n$,
		simply observe that $T = CT^*C$ and $T = -KT^*K$ where
		\begin{equation*}
			C = \minimatrix{0}{J}{J}{0}, \quad
			K = \minimatrix{0}{J}{-J}{0}.
		\end{equation*}
		Now appeal to Lemma \ref{LemmaUECSM} and Lemma \ref{LemmaUESHM}.
	\end{proof}

	\begin{Example}\label{Example6x6}
		The matrix 	
		\begin{equation*}
			T = 
			\left(
			\begin{array}{ccc|ccc}
				0 & 1 & 0 & 0 & 0 & 0 \\
				0 & 0 & 2 & 0 & 0 & 0 \\
				0 & 0 & 0 & 0 & 0 & 0 \\
				\hline
				0 & 0 & 0 & 0 & 0 & 0 \\
				0 & 0 & 0 &1 & 0 & 0 \\
				0 & 0 & 0 & 0 & 2 & 0 \\
			\end{array}
			\right)
		\end{equation*}
		is of the form \eqref{eq-UET2x2} and hence is UET by Lemma \ref{LemmaBlock}.
		However, $T$ has Halmos' matrix \eqref{eq-HalmosExample} as a direct summand
		and that specific matrix is known \emph{not} to be UET \cite[Pr.~159]{Halmos}.  
		This example indicates that we cannot take a block diagonal matrix which is
		UET and conclude that the direct summands are also UET.  
	\end{Example}

	Lemma \ref{LemmaBlock} asserts that matrices of the form $A \oplus A^t$
	are both UECSM and UESHM.  However, by the nature of their construction 
	such matrices are reducible.  On the other hand, 
	for \emph{irreducible} matrices we have the following lemma.	

	\begin{Lemma}\label{LemmaIntersection}
		An irreducible matrix cannot be both UESHM and UECSM.
	\end{Lemma}

	\begin{proof}
		Suppose toward a contradiction that $T$ is
		an irreducible matrix which is both UECSM and UESHM.  Since $T^*$
		also shares these same properties, by Lemma \ref{LemmaUECSM} and 
		Lemma \ref{LemmaUESHM} there is a conjugation $C$ and an anticonjugation $K$ such that
		\begin{equation}\label{eq-CTCKTK}
			CTC = T^* =  -KTK.
		\end{equation}
		Since $C^2 = I$ and $K^2 = -I$ we conclude from \eqref{eq-CTCKTK} that
		$(CK)T = T(CK)$ whence $T$ commutes with the unitary operator $U = CK$.
		Since $T$ is irreducible we conclude that $CK = \alpha I$
		for some unimodular constant $\alpha$.  Multiplying both sides of the preceding by $C$
		we obtain $K  = \overline{\alpha}C$ from which we obtain the contradiction
		\begin{equation*}
			-I = K^2 =  (\overline{\alpha} C)(\overline{\alpha}C)= |\alpha|^2 C^2 = I. \qedhere
		\end{equation*}
	\end{proof}

	We should pause to remark that Lemma \ref{LemmaIntersection}
	does not give a contradiction in the $2 \times 2$ case.  
	Although every $2 \times 2$
	matrix is UECSM by Lemma \ref{Lemma2x2UECSM}, 
	every $2 \times 2$ matrix which is UESHM must actually be scalar
	(and hence reducible) since it is of the form $A \oplus A^t$ for some $1 \times 1$
	matrix $A$.  
	
	The following proposition provides a partial converse to the implication
	\eqref{eq-Implication}:

	\begin{Proposition}\label{PropositionIrreducible}
		If $T \in M_n(\C)$ is irreducible and UET, then precisely one of the following is true:
		\begin{enumerate}\addtolength{\itemsep}{0.5\baselineskip}
			\item $T$ is UECSM,
			\item $T$ is UESHM.
		\end{enumerate}
		\smallskip
		If in addition $n$ is odd, then $T$ must be UECSM.
	\end{Proposition}

	\begin{proof}
		Suppose that $T = UT^tU^*$ for some unitary matrix $U$.
		Taking the transpose of the preceding we obtain
		$T^t = \overline{U}TU^t$ whence $(U\overline{U})T = T(U\overline{U})$.
		Since $T$ is irreducible and $U\overline{U}$ is unitary, it follows that
		$U\overline{U} = \alpha I$ for some $|\alpha|=1$.  From this we conclude
		that $U = \alpha U^t$ whence $U^t = \alpha U$ follows upon transposition.  Putting this all together,
		we conclude that $U = \alpha^2 U$ and consequently $\alpha^2 = 1$. 
		If $\alpha = 1$, then $U = U^t$ whence $T$ is UECSM by Lemma \ref{LemmaUECSM}.
		On the other hand, if $\alpha = -1$, then $U = - U^t$ whence 
		$T$ is UESHM by Lemma \ref{LemmaUESHM}.  The final statement also follows
		from Lemma \ref{LemmaUESHM}.
	\end{proof}		

	Putting the preceding material together we obtain the following.

	\begin{Proposition}\label{PropositionAAt}
		If $T = A \oplus A^t$, then $T$ is unitarily equivalent to a
		direct sum of matrices, each of which is either
		\begin{enumerate}\addtolength{\itemsep}{0.5\baselineskip}
			\item an irreducible complex symmetric matrix,
			\item an irreducible skew-Hamiltonian matrix,
			\item a block matrix of the form $A_i \oplus A_i^t$ where $A_i$ is irreducible
				and neither UECSM nor UESHM.
		\end{enumerate}
	\end{Proposition}

	\begin{proof}
		Suppose that $A$ is irreducible.  
		If $A$ is not UET, then we can conclude that $A$ is neither UECSM (by Lemma \ref{LemmaUECSM})
		nor UESHM (by Lemma \ref{LemmaUESHM}).  Thus $T$ is already 
		of the form (iii).  On the other hand, if $A$ is UET, then $A$ is either UECSM or UESHM
		by Proposition \ref{PropositionIrreducible}.  In either case, $T$ is of the desired form.
		
		If $A$ is reducible, then write $A \cong A_1 \oplus A_2 \oplus \cdots \oplus A_r$
		where $A_1,A_2,\ldots, A_r$ are irreducible and note that
		\begin{align*}
			T  
			&= A \oplus A^t \\
			&\cong (A_1 \oplus A_2 \oplus \cdots \oplus A_r) 
				\oplus (A_1^t \oplus A_2^t \oplus \cdots \oplus A_r^t) \\
			&\cong (A_1 \oplus A_1^t) \oplus (A_2 \oplus A_2^t) \oplus \cdots \oplus (A_r \oplus A_r^t).
		\end{align*}
		Now apply the first portion of the proof to each of the matrices $T_i = A_i \oplus A_i^t$.
	\end{proof}

	\begin{figure}
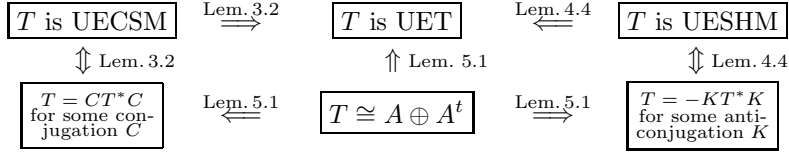

		\begin{equation*}
			\begin{array}{ccccc}
				\boxed{\text{$T$ is UECSM}} &\overset{\text{Lem.\,\ref{LemmaUECSM}}  }{\Longrightarrow}
				& \boxed{\text{$T$ is UET}} &\overset{\text{Lem.\,\ref{LemmaUESHM}}}{\Longleftarrow} 
				&\boxed{\text{$T$ is UESHM}}\\[2pt]
				\,\,\,\qquad \Updownarrow \text{\scriptsize Lem.\,\ref{LemmaUECSM}}& & 
				\quad\qquad \Uparrow \text{\scriptsize Lem.~\ref{LemmaBlock}}
				 & & 
				\,\,\,\,\qquad \Updownarrow\text{\scriptsize Lem.\,\ref{LemmaUESHM}}\\[2pt]
				\boxed{\substack{\text{$T=CT^*C$}\\ \text{for some con-} \\ \text{jugation $C$}} }
				&\overset{\text{Lem.\,\ref{LemmaBlock}}}{\Longleftarrow}& \boxed{ T \cong A \oplus A^t } 
				&\overset{\text{Lem.\,\ref{LemmaBlock}}}{\Longrightarrow}&
				\boxed{ \substack{\text{$T= -KT^*K$}\\ \text{for some anti-} \\ \text{conjugation $K$} } }\\
			\end{array}
		\end{equation*}	
		\caption{\footnotesize Relationship between the classes UET, UECSM, and UESHM.
		In general, the one-way implications cannot be reversed.}
	\end{figure}	
	
\section{SHMs in dimensions $2$, $4$, and $6$ are reducible}\label{SectionReducible}

	It turns out that skew-Hamiltonian matrices in dimensions $2$, $4$, and $6$
	are reducible whereas irreducible SHMs exist 
	in dimensions $8,10,12,\ldots$ (see Section \ref{SectionIrreducible}).
	This issue is highly nontrivial, for it is known that every skew-Hamiltonian matrix is
	\emph{similar} to a block matrix of the form $A \oplus A^t$ \cite{Waterhouse}.
	
	Since every $2 \times 2$ SHM is obviously scalar, we consider in this section
	the $4\times 4$ and $6 \times 6$ cases.

	\begin{Lemma}\label{Lemma4x4}
		If $T \in M_4(\C)$ is SHM, then $T$ is reducible.
	\end{Lemma}
		
	\begin{proof}
		Suppose that $T \in M_4(\C)$ is of the form \eqref{eq-SHM}.
		By interchanging the second and third rows, and then the second and third columns of $T$, 
		we may further assume that $T$ is of the form
		\begin{equation}\label{eq-T4x4UESHM}
			T = \minimatrix{\lambda_1 I}{X}{ \widetilde{X} }{ \lambda_2 I },
		\end{equation}
		where $\lambda_1, \lambda_2 \in \C$, $X$ is a $2 \times 2$ matrix, and 
		$\widetilde{X}$ denotes the \emph{adjugate} of $X$.  In other words, we have
		\begin{equation*}
			X = \minimatrix{a}{b}{c}{d}, \qquad \widetilde{X} = \minimatrix{d}{-b}{-c}{a},
		\end{equation*}
		where $X\widetilde{X} = \widetilde{X}X = (\det X)I$.  Note that the adjugate operation
		satisfies $\widetilde{XY} = \widetilde{Y} \widetilde{X}$ and $\widetilde{X} = (\det X)X^{-1}$
		if $X$ is invertible.
		
		Now write $X = UDW^*$ where $U,W$ are unitary matrices satisfying $\det U = \det W = 1$
		and $D$ is a diagonal matrix having \emph{complex} entries 
		(this factorization can easily  be obtained from the singular
		value decomposition of $X$).  Plugging this into \eqref{eq-T4x4UESHM} we find that
		\begin{align*}
			T
			&= \minimatrix{\lambda_1 I}{UDW^*}{ \widetilde{UDW^*} }{ \lambda_2 I } \\
			&= \minimatrix{\lambda_1 I}{UDW^*}{ W \widetilde{D}U^* }{ \lambda_2 I } \\
			&= \minimatrix{U}{0}{0}{W}\minimatrix{\lambda_1 I}{D}{ \widetilde{D} }{\lambda_2 I}
			\minimatrix{U^*}{0}{0}{W^*} .
		\end{align*}
		Writing
		\begin{equation*}
			D = \minimatrix{ \xi_1}{0}{0}{\xi_2},
		\end{equation*}
		where $\xi_1,\xi_2 \in \C$, we find that 
		\begin{equation*}
			T \cong
			\minimatrix{\lambda_1 I}{D}{ \widetilde{D} }{\lambda_2 I}
			=
			\left(
			\begin{array}{cc|cc}
				\lambda_1 & 0& \xi_1 & 0 \\
				0 & \lambda_1 & 0 & \xi_2 \\
				\hline
				\xi_2 & 0 & \lambda_2 & 0 \\
				0 & \xi_1 & 0 & \lambda_2
			\end{array}
			\right)
			\cong
			\left(
			\begin{array}{cc|cc}
				\lambda_1 & \xi_1 & 0 & 0 \\
				\xi_2 & \lambda_2 & 0 & 0 \\
				\hline
				0 & 0 & \lambda_1 & \xi_2 \\
				0 & 0 & \xi_1 & \lambda_2
			\end{array}
			\right). \qedhere
		\end{equation*}
	\end{proof}
	
	While it is true that every $6 \times 6$ SHM is reducible, the proof is not nearly as simple
	as that of Lemma \ref{Lemma4x4}.	
	Unfortunately, we were unable to come up with a proof that did not involve a significant
	amount of symbolic computation.  Nevertheless, the techniques involved are relatively
	simple and the motivated reader should have no trouble verifying the calculations
	described below.

	\begin{Lemma}\label{Lemma6x6}
		If $T \in M_6(\C)$ is SHM, then $T$ is reducible.
	\end{Lemma}

	\begin{proof}
		Let $T$ be a $6 \times 6$ matrix of the form \eqref{eq-SHM}.
		We intend to show that $T$ commutes with a non-scalar 
		selfadjoint matrix $Q$.  However, attempting to consider
		$QT = TQ$ as a system of $36 \times 2 = 72$ real equations in $15\times2+6 = 36$ 
		real variables is not computationally feasible since the resulting $72 \times 36$ system 
		is symbolic, rather than numeric.  
		Several simplifications are needed before our problem becomes tractable.

		First let us note that if $U \in M_3(\C)$ is unitary then
		\begin{equation}\label{eq-Appeal}
			\minimatrix{U}{0}{0}{\overline{U}} \minimatrix{A}{B}{D}{A^t} \minimatrix{U^*}{0}{0}{U^t}
			= \minimatrix{  UAU^* }{  UBU^t }{  \overline{U}DU^* }{  (UAU^*) ^t }
		\end{equation}
		where $UBU^t$ and $\overline{U}DU^*$ are skew-symmetric.  Without loss
		of generality, we may therefore assume that $A$ is upper-triangular.  Furthermore,
		by scaling and subtracting a multiple of the identity from $T$ we can also assume 
		that at most one of the diagonal entries of $A$ is non-real while the others are either $0$ or $1$.

		Appealing to \eqref{eq-Appeal} again where $U$ is now a suitable diagonal unitary matrix, we can further
		arrange things so that that the skew-symmetric matrix $UBU^t$ has only real entries.  
		The nonzero off-diagonal entries of the upper-triangular matrix $UAU^*$ may 
		change depending on $U$, but the diagonal entries are unaffected.
		Thus we may further assume that $B$ has only real entries.

		Now recall from \eqref{eq-Symplectic} that $T$ is SHM if and only if $T = -\Omega T^t \Omega$
		where $\Omega$ denotes the matrix \eqref{eq-Omega}.  If $Q = Q^*$ and $QT = TQ$, then
		clearly $Q^t T^t = T^t Q^t$ whence $\Omega Q^t \Omega$ also
		commutes with $T$.  Let us therefore consider selfadjoint matrices $Q$ which satisfy
		$Q = \Omega Q^t \Omega$.  In other words, we restrict our attention to matrices of the form
		\begin{equation}\label{eq-Commuter}
			Q = \minimatrix{X}{Y}{Y^*}{-X^t}
		\end{equation}
		where $X = X^*$ and $Y = Y^t$.  There are now a total of $9 + 12 = 21$ real
		unknowns arising from the components of $X$ and $Y$.

		Expanding out the system $QT = TQ$ we obtain 
		\begin{align}
			XA-AX+YD-BY^* &=0, \label{eq-6x6-01} \\
			XB+BX^t + YA^t - AY &=0, \label{eq-6x6-02} \\
			DX + X^t D + A^t Y^* - Y^*A &=0, \label{eq-6x6-03} \\
			A^t X^t - X^t A^t -DY + Y^*B &=0, \label{eq-6x6-04}
		\end{align}
		which yields a system of $72$ real equations in $21$ unknowns.
		However it is clear that \eqref{eq-6x6-01} and \eqref{eq-6x6-04} are transposes
		of each other and hence \eqref{eq-6x6-04} can be ignored.  Moreover,
		\eqref{eq-6x6-02} and \eqref{eq-6x6-03} are skew-symmetric and hence we obtain a system
		of $18 + 6 + 6 = 30$ real equations in $21$ unknowns.  

		With the preceding reductions in hand, it becomes possible to compute the rank of the system symbolically
		via \texttt{Mathematica}.  In particular, the \texttt{MatrixRank} command computes the
		rank of a symbolic matrix under the assumption that the distinct symbols appearing 
		as coefficients in the system are linearly independent.
		By considering separately the cases where $A$ has either one, two, or three
		distinct eigenvalues one can conclude that the rank of our system
		is always $\leq 20$ whence a nontrivial solution to our problem exists.
		By \eqref{eq-Commuter} it is clear 
		that the resulting $Q$ cannot be a multiple of the identity whence $T$ is reducible.
	\end{proof}

\section{Irreducible SHMs exist in dimensions $n=8,10,12,\ldots$}\label{SectionIrreducible}

	We now turn our attention to the task of proving that irreducible SHMs exist in dimensions $n=8,10,12,\ldots$
	As we have mentioned before, this is nontrivial due to the fact that every skew-Hamiltonian
	matrix \emph{similar} to a block matrix of the form $A \oplus A^t$ \cite{Waterhouse}.
	
	\begin{Proposition}\label{Proposition8x8}
		For each $d \geq 4$ the $2d \times 2d$ matrix 
		\begin{equation*}
			T = \minimatrix{A}{B}{0}{A},
		\end{equation*}
		where
		\begin{equation}\label{eq-AB}
			A = 
			\begin{pmatrix}
				1 & & &\\
				 & 2 & & \\[-3pt]
				 & & \ddots & \\[-3pt]
				 & & & d
			\end{pmatrix}, 
			\qquad 
			B=
			\begin{pmatrix}
				0 & 1 & 1 &\cdots & 1 &  1\\[-3pt]
				-1 & 0 & 1& \ddots &  1 &  1\\[-3pt]
				-1 & -1 & 0 & \ddots & 1 &  1\\[-3pt]
				\vdots & \ddots& \ddots & \ddots &  \ddots &  \vdots\\[-3pt]
				-1 & -1 &  -1 & \ddots & 0& 1 \\
				-1 & -1 &  -1 & \cdots & -1& 0
			\end{pmatrix},
		\end{equation}
		is skew-Hamiltonian and irreducible.
	\end{Proposition}

	\begin{proof}
		Since $A = A^t$ and $B^t = -B$, it is clear that $T$ is SHM.
		It therefore suffices to prove that $T$ is irreducible.  To this end, suppose that 
		$Q$ is a selfadjoint matrix which satisfies $QT = TQ$.  Writing
		\begin{equation} \label{eq-QXXYZ}
			Q= \minimatrix{X}{Y}{Y^*}{Z}, 
		\end{equation}
		where $X = X^*$ and $Z = Z^*$, we find that
		\begin{align}
			XA &= AX+BY^*, \label{eq-8x8-01}\\
			XB+YA &=AY+BZ,\label{eq-8x8-02}\\
			Y^*A &= AY^*,\label{eq-8x8-03}\\
			Y^*B + ZA &=AZ.\label{eq-8x8-04}
		\end{align}
		Examining \eqref{eq-8x8-03} entry-by-entry reveals that $Y^*$ is diagonal.
		In particular, $YA = AY$ and hence
		\begin{equation}\label{eq-XBBZ}
			XB = BZ
		\end{equation}
		follows from \eqref{eq-8x8-02}.
		By \eqref{eq-8x8-01} we have 
		\begin{equation}\label{eq-XAAXBY*}
			\underbrace{XA - AX}_{\text{skew-selfadjoint}} = BY^*.
		\end{equation}
		Since $Y^*$ is diagonal, a short computation
		shows that $BY^*$ is skew-selfadjoint if and only if
		$Y^* = \alpha I$ for some $\alpha \in \R$ (this requires $d \geq 3$).
		We may therefore rewrite \eqref{eq-XAAXBY*} as
		\begin{equation}\label{eq-XAAX}
			XA - AX = \alpha B.
		\end{equation}
		Equation \eqref{eq-8x8-04} now assumes the similar form
		\begin{equation}\label{eq-ZAAZ}
			ZA - AZ = -\alpha B.
		\end{equation}
		Adding \eqref{eq-XAAX} and \eqref{eq-ZAAZ} together we find that
		\begin{equation*}
			(X+Z)A  = A(X+Z)
		\end{equation*}
		whence, since $A$ is diagonal and has distinct eigenvalues, the matrix
		$X+Z = D$ is also diagonal.  Plugging this into \eqref{eq-XBBZ} yields
		\begin{equation*}
			\underbrace{XB +BX}_{\text{skew-selfadjoint}} = BD.
		\end{equation*}
		The same reasoning employed in analyzing \eqref{eq-XAAXBY*} now reveals that $D = 2\delta I$
		for some $\delta \in \R$.  Since $X +Z = 2\delta I$ we conclude that
		\begin{equation}\label{eq-Subtract}
			(X -  \delta I) =- (Z -  \delta I).
		\end{equation}
		At this point we observe that $Q+  \delta I$ also commutes with $T$,
		and hence upon making the substitutions  $X \mapsto X + \delta I$ and $Z \mapsto Z + \delta I$ 
		in \eqref{eq-QXXYZ} we may assume that $X = -Z$.  Plugging this into \eqref{eq-XBBZ}
		we see that
		\begin{equation}\label{eq-XBBX}
			XB = -BX.
		\end{equation}
		From equations \eqref{eq-XAAX} and \eqref{eq-XBBX} we shall derive a number of
		constraints upon the entries of $X$ which can be shown to be mutually incompatible 
		unless $X = Y = 0$.

		Examining \eqref{eq-XAAX} entry-by-entry, we find that the $ij$th entry $x_{ij}$ of $X$ is given by
		\begin{equation*}
		x_{ij} = \frac{ \alpha }{|j - i|}
		\end{equation*}
		for $i \neq j$.  Thus 
		\begin{equation}\label{eq-Xplicit}
			X = 
			\begin{pmatrix}
				x_{11} & \alpha & \frac{\alpha}{2} & \frac{ \alpha}{3} & \cdots & \frac{\alpha}{d-1} \\
				\alpha & x_{22} & \alpha & \frac{\alpha}{2} & \cdots & \frac{ \alpha}{d-2} \\
				\frac{\alpha}{2} & \alpha & x_{33} & \alpha & \cdots & \frac{ \alpha }{d-3} \\
				\frac{\alpha}{3} &\frac{\alpha}{2}  & \alpha & x_{44} &  \cdots & \frac{ \alpha }{d-4} \\
				\vdots & \ddots & \ddots & \ddots & \ddots & \vdots \\
				\frac{ \alpha}{d-1} & \frac{ \alpha}{d-2}& \frac{ \alpha}{ d-3} & \frac{ \alpha}{d-4} & \cdots & x_{dd}
			\end{pmatrix}
		\end{equation}
		where the diagonal entries $x_{11}, x_{22}, \ldots, x_{dd}$ are to be determined.

		On the other hand, from \eqref{eq-XBBX} it follows that
		\begin{equation}\label{eq-UXB}
			(UX)B = B(UX)
		\end{equation}
		since $UBU= -B$ where
		\begin{equation*}
			U = 
			\begin{pmatrix}
				&&1\\[-3pt]
				&\iddots&\\[-3pt]
				1&&
			\end{pmatrix}.
		\end{equation*}
		From \eqref{eq-UXB} we wish to conclude that $UX = p(B)$ for some polynomial $p(z)$.
		This will follow if we can show that $B$ has distinct eigenvalues.

		A short computation first reveals 
		\begin{equation*}
			2(B+I)^{-1} = 
			\begin{pmatrix}
				1 & -1 & 0 & 0 & 0 \\
				 0 & 1 & -1 & 0  &0\\
				 0 & 0 & 1 & -1 &0 \\[-5pt]
				 \vdots & \vdots & \vdots & \ddots & \vdots \\
				0 & 0 & 0 & 1 &-1\\
				 1 & 0 & 0 & 0 &1\\
			\end{pmatrix}
		\end{equation*}
		whence
		\begin{equation*}
			I - 2(B+I)^{-1} = 
			\underbrace{
			\begin{pmatrix}
				0 & 1 & 0 & 0 & 0 \\
				 0 & 0 & 1 & 0  &0\\
				 0 & 0 & 0 & 1 &0 \\[-5pt]
				 \vdots & \vdots & \vdots & \ddots & \vdots \\
				0 & 0 & 0 & 0 &1\\
				 -1 & 0 & 0 & 0 &0\\
			\end{pmatrix}
			}_L.
		\end{equation*}
		Cofactor expansion along the first column shows that 
		$L$ has the characteristic polynomial
		$z^d +1$, whose roots are precisely the $d$th roots of $-1$.  From this we conclude
		that the skew-symmetric matrix $B$ has $d$ distinct eigenvalues (namely the numbers
		$(1+\zeta_i)/(1 - \zeta_i)$ where $\zeta_1,\zeta_2,\ldots,\zeta_d$ are the $d$th roots of $-1$).
		
		Since $1$ is not an eigenvalue of $L$ it follows that
		\begin{equation*}
			B = 2(I - L)^{-1} - I = q(L)
		\end{equation*}
		is a polynomial in $L$ (the fact that $(I - L)^{-1}$ is a polynomial in $L$ follows
		from the Cayley-Hamilton theorem).  Thus $UX = p(B) = p(q(L))$ is also a polynomial in $L$.
		Now observe that $L^d = - I$ and that each of the matrices $I,L,L^2,\ldots,L^{d-1}$
		is a Toeplitz matrix.  Therefore $UX$ is a Toeplitz matrix
		whence $X$ is a Hankel matrix.
		
		Recalling that $d \geq 4$ and looking back to the explicit formula \eqref{eq-Xplicit}
		for $X$, we see that $X$ cannot be a Hankel matrix unless 
		$\alpha = 0$ and $x_{11} = x_{22} = \cdots = x_{dd} = 0$.
		In other words, it must be the case that $X= Y = 0$ whence $Q = 0$
		(note that we considered $Q + \delta I$ in place of $Q$
		after the step \eqref{eq-Subtract}).  Thus $T$ is irreducible.
	\end{proof}

	Observe that the argument above fails if $d \leq 3$ since setting $x_{22} = \frac{\alpha}{2}$
	in \eqref{eq-Xplicit} yields a Hankel matrix.  In this case, one can verify directly 
	that $T$ commutes with 
	\begin{equation*}
		\left(
		\begin{array}{ccc|ccc}
			-1 & 1 & \frac{1}{2} & 1 & 0 & 0 \\
			1 & \frac{1}{2} & 1 & 0 & 1 & 0\\
			\frac{1}{2}& 1 & -1 & 0 & 0 & 1 \\
			\hline
			1 &0 & 0 & 1 & -1 & -\frac{1}{2} \\
			0 & 1 & 0 & -1 & -\frac{1}{2} & -1 \\
			0 & 0 & 1 & -\frac{1}{2} & -1 & 1
		\end{array}
		\right).
	\end{equation*}
	As anticipated in the proof of Lemma \ref{Lemma6x6}, the preceding
	matrix is of the form \eqref{eq-Commuter}.

\section{Proof of Theorem \ref{TheoremMain}}\label{SectionProof01}

This entire section is devoted to the proof of Theorem \ref{TheoremMain}.
For the sake of readability, it is divided into several subsections.

\subsection{The simple direction}
	One implication of the proof is now trivial.  If $T$ is unitarily equivalent to 
	a direct sum of CSMs, SHMs, and matrices of the form $A \oplus A^t$,
	then $T$ is UET by Lemmas \ref{LemmaUECSM}, \ref{LemmaUESHM},
	and \ref{LemmaBlock}.

\subsection{Initial setup}
	We now focus on the more difficult implication of Theorem \ref{TheoremMain}.
	If $T \in M_n(\C)$ is UET, then there exists a unitary matrix $U$ such that
	\begin{equation}\label{eq-Original}
		T = UT^tU^*.
	\end{equation}
	Taking the transpose of the preceding and solving for $T^t$ we obtain
	$T^t = \overline{U}TU^t$ whence
	\begin{equation}\label{eq-Commutant}
		(U\overline{U})T = T(U\overline{U}).
	\end{equation}
	In light of \eqref{eq-Commutant} we are therefore led to consider the unitary matrix 
	\begin{equation*}
		V = U\overline{U}.
	\end{equation*}
	The following lemma lists several restrictions on the eigenvalues of $V$ which 
	will be useful in what follows.  We denote by $\sigma(A)$ the set of eigenvalues of a matrix $A \in M_n(\C)$.
	
	\begin{Lemma}\label{LemmaV}
		If $U$ is a unitary matrix and $V = U \overline{U}$, then $\det V= 1$
		and $\sigma(V) = \overline{ \sigma(V)}$.
		In particular, the eigenvalues of $V$ are restricted to:
		\begin{enumerate}\addtolength{\itemsep}{0.5\baselineskip}
			\item $1$,
			\item $-1$, with even multiplicity,
			\item complex conjugate pairs $\lambda,\overline{\lambda}$ where $\lambda\neq \pm 1$ 
					and both $\lambda$ and $\overline{\lambda}$ have the same multiplicity.
		\end{enumerate}
	\end{Lemma}

	\begin{proof}
		Since $U$ is unitary it follows that
		\begin{equation*}
			\det V = \det U\overline{U}  = (\det U)(\det \overline{U})
			= (\det U) \overline{(\det U)} = |\det U|^2 = 1.
		\end{equation*}
		Since $U$ is invertible, $U\overline{U}$ and $\overline{U}U$ have the same characteristic polynomial whence
		\begin{equation*}
			\dim \ker ( U \overline{U} - z I)
			= \dim \ker(\overline{U}U - \overline{z} I) 
			= \dim \ker(U\overline{U} - \overline{z} I)
		\end{equation*}
		holds for all $z \in \C$.  Among other things, this establishes (iii) and shows that $\sigma(V) = \overline{ \sigma(V)}$.
		Conditions (i) and (ii) now follow from (iii) and the fact that $\det V = 1$.
	\end{proof}

	Now observe that 
	\begin{equation}\label{eq-UVVU}
		U \overline{V} = V U
	\end{equation}
	holds since both sides of \eqref{eq-UVVU} equal
	$U \overline{U}U$.  By the Spectral Theorem, there exists
	a unitary matrix $W$ so that 
	\begin{equation}\label{eq-VMatrix}
		V = W^*DW
	\end{equation}
	where (by Lemma \ref{LemmaV}) $D$ is a block diagonal matrix of the form
	\begin{equation}\label{eq-DMatrix}
		D = 
		\left(
		\begin{array}{c|c|cc|c|cc}
		I &&&&&&\\
		\hline
		 &-I&&&&&\\
		 \hline
		 &&\lambda_1 I&&&&\\
		& &&\overline{\lambda_1} I&&&\\
		\hline
		&&&&\ddots&& \\
		\hline
		&&&&&\lambda_r I& \\
		&&&&&&\overline{\lambda_r} I\\
		\end{array}
		\right)
	\end{equation}
	for some distinct unimodular constants $\lambda_1,\lambda_2,\ldots,\lambda_r$ such that
	\begin{enumerate}\addtolength{\itemsep}{0.5\baselineskip}
		\item $\lambda_i \neq \pm 1$ for $i = 1,2,\ldots,r$,
		\item $\lambda_i \neq \overline{\lambda_j}$ for $1 \leq i,j \leq r$.
	\end{enumerate}
	Let us make several remarks about the matrix
	\eqref{eq-DMatrix}.  First, it is possible that some of the blocks may be absent depending
	upon $V$.
	Second, the blocks corresponding to conjugate eigenvalues $\lambda_i$ and $\overline{\lambda_i}$
	must be of the same size by Lemma \ref{LemmaV}.
	
\subsection{The matrix $Q$}	

	Substituting \eqref{eq-VMatrix} into \eqref{eq-UVVU} we find that
	\begin{equation}\label{eq-QDDQ}
		Q\overline{D} = D Q
	\end{equation}
	where
	\begin{equation}\label{eq-QMatrix}
		Q = WUW^t
	\end{equation}
	is unitary.
	In light of \eqref{eq-QDDQ} and the structure of $D$ given in \eqref{eq-DMatrix},
	a short matrix computation reveals that
	\begin{equation}\label{eq-QMatrixExplicit}
		Q = 
		\left(
		\begin{array}{c|c|cc|c|cc}
		Q_+ &&&&&&\\
		\hline
		 &Q_-&&&&&\\
		 \hline
		 &&0&Y_1&&&\\
		& &X_1&0&&&\\
		\hline
		&&&&\ddots&& \\
		\hline
		&&&&&0&Y_r \\
		&&&&&X_r&0\\
		\end{array}
		\right)
	\end{equation}
	where $Q_+,Q_-, X_1,Y_1,\ldots, X_r,Y_r$ are unitary matrices.
	Next observe that
	\begin{equation*}
		Q\overline{Q} = (WUW^t)(\overline{W} \overline{U}W^*) = WU\overline{U}W^* = WVW^* = D
	\end{equation*}
	by \eqref{eq-QMatrix} and \eqref{eq-VMatrix}.  Using the fact that $Q$ is unitary we conclude 
	from the preceding that
	\begin{equation}\label{eq-QDQt}
		Q = DQ^t.
	\end{equation}
	This gives us further insight into the structure of $Q$.  Using the block matrix
	decompositions \eqref{eq-DMatrix} and \eqref{eq-QMatrixExplicit}
	and examining \eqref{eq-QDQt} block-by-block we conclude that
	\begin{enumerate}\addtolength{\itemsep}{0.5\baselineskip}
		\item $Q_+ = Q_+^t$ (i.e., $Q_+$ is complex symmetric and unitary),
		\item $Q_- = -Q_-^t$ (i.e., $Q_-$ is skewsymmetric and unitary),
		\item $Y_i = \lambda_i X_i^t$ where $X_i$ is unitary and $\lambda_i \neq \pm 1$
			for $i = 1,2,\ldots,r$.
	\end{enumerate}
	Thus we may write
	\begin{equation}\label{eq-QMatrixExplicit2}
		Q = 
		\left(
		\begin{array}{c|c|cc|c|cc}
		Q_+ &&&&&&\\
		\hline
		 &Q_-&&&&&\\
		 \hline
		 &&0&\lambda_1 X_1^t&&&\\
		& &X_1&0&&&\\
		\hline
		&&&&\ddots&& \\
		\hline
		&&&&&0&\lambda_r X_r^t \\
		&&&&&X_r&0\\
		\end{array}
		\right).
	\end{equation}
	It is important to note that some of the blocks in \eqref{eq-QMatrixExplicit2} may be absent,
	depending upon the decomposition \eqref{eq-DMatrix} of $D$.  

\subsection{Simplifying the problem}\label{SectionSimplify}
	
	Now that we have a concrete description of $Q$ in hand, let us return to our
	original equation \eqref{eq-Original}.  In light of \eqref{eq-QMatrix} we see that
	\begin{equation*}
		U = W^*Q\overline{W}
	\end{equation*}
	whence by \eqref{eq-Original} it follows that
	\begin{equation*}
		T = (W^*Q\overline{W})T^t (W^t Q^*W) = W^* Q (WTW^*)^t  Q^*W.
	\end{equation*}
	Rearranging this reveals that
	\begin{equation*}
		(WTW^*) Q
		= Q(WTW^*)^t .
	\end{equation*}
	Since $T \cong WTW^*$, it follows that in order to characterize, 
	up to unitary equivalence, those matrices $T$ which are UET
	we need only consider those $T$ which satisfy
	\begin{equation}\label{eq-New}
		TQ =QT^t
	\end{equation}
	where $Q$ is a unitary matrix of the special form \eqref{eq-QMatrixExplicit2}	
	satisfying conditions (i), (ii), and (iii) above.	
	
\subsection{Describing $T$}
	
	Taking the transpose of \eqref{eq-New}, solving for $T^t$, and substituting the result back into
	\eqref{eq-New} reveals that $(Q\overline{Q})T = T(Q \overline{Q})$.
	Thus $T$ is actually block-diagonal, the sizes of the corresponding blocks
	being determined by the decomposition \eqref{eq-QMatrixExplicit2} of $Q$ itself.
	Returning to \eqref{eq-New},	
	we find that $T$ has the form	
	\begin{equation}\label{eq-TMatrix}
		T = 
		\left(
		\begin{array}{c|c|cc|c|cc}
		T_+ &&&&&&\\
		\hline
		 &T_-&&&&&\\
		 \hline
		 &&A_1&0&&&\\
		&&0&X_1A_1^tX_1^*&&&\\
		\hline
		&&&&\ddots&& \\
		\hline
		&&&&&A_r&0 \\
		&&&&&0&X_r A_r^t X_r^*\\
		\end{array}
		\right)
	\end{equation}
	where
	\begin{enumerate}\addtolength{\itemsep}{0.5\baselineskip}
		\item $T_+ = Q_+ T_+^t Q_+^*$.  Since $Q_+$ is a symmetric unitary matrix,
			Lemma \ref{LemmaUECSM} implies that $T_+$ is UECSM.

		\item $T_- = Q_- T_-^t Q_-^*$.  Since $Q_-$ is a skewsymmetric unitary matrix,
			Lemma \ref{LemmaUESHM} implies that $T_-$ is UESHM.

		\item $A_1,A_2,\ldots,A_r$ are arbitrary.
	\end{enumerate}
	Obtaining (i) and (ii) from \eqref{eq-New} is straightforward.  Let us say a few words about (iii).  Examining
	\eqref{eq-New} we obtain $r$ equations of the form
	\begin{equation*}
		\minimatrix{A_i}{B_i}{C_i}{D_i} \minimatrix{0}{\lambda_i X_i^t}{X_i}{0} = 
		\minimatrix{0}{\lambda_i X_i^t}{X_i}{0}\minimatrix{A_i^t}{C_i^t}{B_i^t}{D_i^t}.
	\end{equation*}
	This leads to $B_i X_i = \lambda_i X_i^t B_i ^t = \lambda_i(B_i X_i)^t$ whence $B_i X_i = \lambda_i^2 B_i X_i$.
	Since $\lambda_i \neq \pm 1$ and $X_i$ is invertible we conclude that $B_i = 0$.  Similarly we find that $C_i = 0$
	and $D_i = X_i A_i^t X_i^*$.
	\smallskip

	It is evident at this point that $T$ is unitarily equivalent to 
	the direct sum of a CSM, an SHM, and blocks of the form
	$A \oplus A^t$.  As usual, some of the blocks in \eqref{eq-TMatrix} may be absent.  The exact
	form of \eqref{eq-TMatrix} depends upon the decomposition \eqref{eq-QMatrixExplicit2}.

	By Proposition \ref{PropositionUECSM}, the block $T_+$ is unitarily equivalent to a
	direct sum of matrices of type I or III.  Similarly, Proposition \ref{PropositionUESHM}
	asserts that $T_-$ is unitarily equivalent to a direct sum of matrices of type II or III.
	Finally, Proposition \ref{PropositionAAt} permits us to decompose the resulting type III
	blocks and the $2 \times 2$ block matrices from \eqref{eq-TMatrix} into
	a direct sum of matrices of type I, II, or III.  
	
	The restrictions on the dimensions of type II and type III blocks follow
	immediately from the results of Sections \ref{SectionReducible}
	and \ref{SectionIrreducible} and Lemma \ref{Lemma2x2UECSM}.	
	This establishes the existence of the desired
	decomposition. 

	For the final statement of Theorem \ref{TheoremMain}, first note that 
	a matrix of type III is reducible and hence cannot belong to the unitary orbit of a type I or type II
	matrix.  That the unitary orbits of a type I matrix and a type II matrix are disjoint
	follows from Lemma \ref{LemmaIntersection}.  \qed

\end{document}